\pdfoutput=1
\RequirePackage{ifpdf}
\ifpdf 
\documentclass[pdftex]{sigma}
\else
\documentclass{sigma}
\fi

\numberwithin{equation}{section}

\begin{document}

\allowdisplaybreaks

\renewcommand{\thefootnote}{$\star$}

\newcommand{\arXivNumber}{1602.02724}

\renewcommand{\PaperNumber}{048}

\FirstPageHeading

\ShortArticleName{Hypergeometric Orthogonal Polynomials with respect to Newtonian Bases}

\ArticleName{Hypergeometric Orthogonal Polynomials\\ with respect to Newtonian Bases\footnote{This paper is a~contribution to the Special Issue
on Orthogonal Polynomials, Special Functions and Applications.
The full collection is available at \href{http://www.emis.de/journals/SIGMA/OPSFA2015.html}{http://www.emis.de/journals/SIGMA/OPSFA2015.html}}}

\Author{Luc VINET~$^\dag$ and Alexei ZHEDANOV~$^\ddag$}

\AuthorNameForHeading{L.~Vinet and A.~Zhedanov}

\Address{$^\dag$~Centre de recherches math\'ematiques, Universit\'e de Montr\'eal,\\
\hphantom{$^\dag$}~P.O.~Box 6128, Centre-ville Station,
Montr\'eal (Qu\'ebec), H3C 3J7 Canada}
\EmailD{\href{mailto:vinet@crm.umontreal.ca}{vinet@crm.umontreal.ca}}

\Address{$^\ddag$~Institute for Physics and Technology, 83114 Donetsk, Ukraine}
\EmailD{\href{mailto:zhedanov@yahoo.com}{zhedanov@yahoo.com}}

\ArticleDates{Received February 08, 2016, in f\/inal form May 07, 2016; Published online May 14, 2016}

\Abstract{We introduce the notion of ``hypergeometric''
 polynomials with respect to Newtonian bases. These polynomials are eigenfunctions ($L P_n(x) = \lambda_n P_n(x)$) of some abstract operator~$L$ which is 2-diagonal in the Newtonian basis $\varphi_n(x)$: $L \varphi_n(x) = \lambda_n \varphi_n(x) + \tau_n(x) \varphi_{n-1}(x)$ with some coef\/f\/icients~$\lambda_n$,~$\tau_n$. We f\/ind the necessary and suf\/f\/icient conditions for the polynomials~$P_n(x)$ to be orthogonal. For the special cases where the sets $\lambda_n$ correspond to the classical grids, we f\/ind the complete solution to these conditions and observe that it leads to the most general Askey--Wilson polynomials and their special and degenerate classes.}


\Classification{42C05; 42C15}

\renewcommand{\thefootnote}{\arabic{footnote}}
\setcounter{footnote}{0}

\section{Introduction}

It is well known that all ``classical'' orthogonal polynomials $P_n(x)$ from the Askey tableau~\cite{KLS} have rather simple expressions of the form
\begin{gather}
P_n(x) = \sum_{s=0}^n W_{ns} \varphi_s(x),
\label{P-hyp-phi}
\end{gather}
where $\varphi_n(x)$ are the Newtonian basis (interpolated) polynomials def\/ined as
\begin{gather*}
\varphi_0 =1, \qquad \varphi_n(x) = (x-a_0) (x-a_1) \cdots (x-a_{n-1}), \qquad n=1,2,\dots 
\end{gather*}
with the real numbers $a_i$ as the interpolation nodes $i=0,1,2,\dots$.

The expansion coef\/f\/icients $W_{ns}$ in these formulas satisfy the two-term recurrence relation
\begin{gather}
W_{n,s+1} = R_{ns} W_{ns}, \label{rec_WR}
\end{gather}
where $R_{ns}$ are simple rational functions either of $n$ or of $q^n$. This property allows to obtain the coef\/f\/icients~$W_{ns}$ explicitly and this leads to the hypergeometric expressions for the classical orthogonal polynomials.

The recurrence relation~\eqref{rec_WR} follows from the fact that the classical orthogonal polyno\-mials $P_n(x)$ satisfy the eigenvalue equation{\samepage
\begin{gather}
L P_n(x) = \lambda_n P_n(x), \label{LPP}
\end{gather}
where $L$ is either a dif\/ferential or dif\/ference operator of second order.}

It turns out that in all these cases the operator $L$ acts on the Newtonian basis according to
\begin{gather}
L \varphi_n(x) = \lambda_n \varphi_n(x) + \tau_n \varphi_{n-1}(x), \label{hyp-2-diag}
\end{gather}
where $\lambda_n$ are the same eigenvalues as in \eqref{LPP} and $\tau_n$ are some coef\/f\/icients.

From \eqref{LPP} and \eqref{hyp-2-diag} follows the recurrence relation \eqref{rec_WR} for these coef\/f\/icients which takes the form
\begin{gather*}
\frac{W_{n,s+1}}{W_{ns}} = \frac{\lambda_n-\lambda_s}{\tau_{s+1}}, 
\end{gather*}
whence
\begin{gather} W_{ns}= W_{n,0}
\frac{(\lambda_n-\lambda_0)(\lambda_n-\lambda_1)\cdots (\lambda_n
-\lambda_{s-1}) }{\tau_1 \tau_2 \cdots \tau_s}, \qquad s=1,2,\dots, n.
\label{expr_W}
\end{gather}
The coef\/f\/icient $W_{n0}$ can be chosen
arbitrarily. One possible choice is $W_{n0}=1$ for all~$n$. This
corresponds to the ``hypergeometric-like'' form of the polynomial~$P_n(x)$. Another choice is the monic form
$P_n(x)=x^n + O(x^{n-1})$. In this case
\begin{gather*} W_{n0} = \frac{\tau_1
\tau_2 \cdots \tau_n}{(\lambda_n-\lambda_0)(\lambda_n-\lambda_1)\cdots
(\lambda_n -\lambda_{n-1}) }. 
\end{gather*} It is then
convenient to present the expansion coef\/f\/icients as follows
\begin{gather}
W_{n,n-k} = \frac{\tau_n \tau_{n-1} \cdots \tau_{n-k+1}}{(\lambda_n-
\lambda_{n-1})(\lambda_n- \lambda_{n-2}) \cdots (\lambda_n-
\lambda_{n-k})}, \qquad k=1,2,\dots, n, \label{W_monic}
\end{gather}
and we thus have the following expression for the polynomial $P_n(x)$:
\begin{gather*}
P_n(x) = \varphi_n(x) + \sum_{k=1}^n {\frac{\tau_n \tau_{n-1} \cdots \tau_{n-k+1}}{(\lambda_n-
\lambda_{n-1})(\lambda_n- \lambda_{n-2}) \cdots (\lambda_n-
\lambda_{n-k})} \varphi_{n-k}(x)}. 
\end{gather*}
One can propose an ``inverse'' problem: assume that the operator $L$ is given in an abstract form by its action~\eqref{hyp-2-diag} with unknown coef\/f\/icients $\lambda_n$, $\tau_n$. The interpolation points~$a_n$ are unknown as well.

The only restrictions for these coef\/f\/icients are:
\begin{enumerate}\itemsep=0pt
\item[(i)] the spectrum is nondegenerate $\lambda_n \ne \lambda_m$ for $n \ne m$;
\item[(ii)] the initial conditions for the coef\/f\/icients $\lambda_n$, $\tau_n$ are
\begin{gather*}
\lambda_0 =\tau_0=a_0=0. 
\end{gather*}
\item[(iii)] In the inf\/inite case $\tau_n \ne 0$ for all $n=1,2,\dots$.
\end{enumerate}

 In the f\/inite case $\tau_n \ne 0$ for all $n=1,2,\dots,N$ and $\tau_{N+1}=0$.

 (For the f\/inite case see the discussion in the next section.)

These conditions are quite natural. Indeed, $\tau_0=0$ follows from the truncation condition of the action of the operator~$L$ on the constant~$\varphi_0(x)=1$. Condition $\lambda_0$ can always be achieved by the appropriate addition of a constant to the operator~$L$. Moreover, we will assume that $a_0=0$. Indeed, if $a_0 \ne 0$, we only need to choose $\tilde \varphi_n(x) = \varphi_n(x-a_0)$; the Newtonian polynomials $\tilde \varphi_n(x)$ have the same properties with respect to the operator~$L$ and the polynomials $P_n(x-a_0)$ and hence we can always assume that $a_0=0$.

Note that the operator $L$ def\/ined by the abstract relation~\eqref{hyp-2-diag} can be considered in the context of the so-called umbral calculus~\cite{Roman} in which case concrete operators (like the derivative operator $\partial_x$) are replaced with their abstract symbols when acting on specif\/ic bases (For instance the abstract umbral derivative operator $\cal D$ can be def\/ined by its symbols $\mu_n$ from the abstract relation ${\cal D} x^n = \mu_n x^{n-1}$; for the ordinary derivative operator one has $\mu_n = n$).

Then it is obvious that the eigenvalue equation~\eqref{LPP} generates a unique system of monic polynomials $P_n(x)$ having explicit expansion coef\/f\/icients \eqref{W_monic}. However, in general, for an arbitrary choice of the coef\/f\/icients $\lambda_n$, $\tau_n$ and $a_n$, the polynomials~$P_n(x)$ will not be orthogonal.

The main problem is to f\/ind the coef\/f\/icients $\lambda_n$, $\tau_n$, $a_n$ such that the eigenpolynomials $P_n(x)$ are orthogonal with respect to a~nondegenerate linear functional $\sigma$
\begin{gather*}
\langle \sigma, P_n(x) P_m(x) \rangle = 0, \qquad n \ne m,
\end{gather*}
where $\sigma$ can be determined in terms of the moments
\begin{gather*}
\langle \sigma, x^n \rangle = c_n, \qquad n=0,1,2,\dots. 
\end{gather*}
The linear functional $\sigma$ is called nondegenerate if the conditions $\Delta_n \ne 0$, $n=0,1,2,\dots$ hold,
where $\Delta_n = |c_{i+k}|_{i,k=0}^n$ are the Hankel determinants constructed from the moments.

Equivalently, the polynomials $P_n(x)$ are orthogonal if and only if they satisfy the three-term recurrence relation
\begin{gather}
P_{n+1}(x) + b_n P_n(x) + u_{n+1} P_{n-1}(x) = xP_n(x), \qquad n=1,2,\dots \label{rec_P}
\end{gather}
with some coef\/f\/icients $b_n$, $u_n$.

The linear functional $\sigma$ is nondegenerate if and only if $u_n \ne 0$, $n=1,2,\dots$~\cite{Chi}. In a~spe\-cial situation when the coef\/f\/icients~$b_n$ are real and $u_n>0$, $n=1,2,\dots$ the polynomials are orthogonal with respect to a positive measure~$d \mu(x)$ on the real axis:
\begin{gather}
\int_{a}^b P_n(x) P_m(x) d \mu(x) = h_n \delta_{nm}, \label{ort_mu}
\end{gather}
where $h_0=1$ and $h_n=u_1 u_2 \cdots u_n$, $n=1,2,\dots$. The integration limits~$a$,~$b$ in~\eqref{ort_mu} may be either f\/inite or inf\/inite.

The simplest case of the monomial basis $\varphi_n(x) = x^n$ was considered in~\cite{Zhe-hyp}. In this instance, a~full classif\/ication was found. This scheme describes all admissible operators~$L$ and corresponding orthogonal polynomials~$P_n(x)$. In a nutshell, the analysis performed in~\cite{Zhe-hyp} led to the following list of of possibilities:
\begin{enumerate}\itemsep=0pt
\item[(i)] $P_n(x)$ are the little $q$-Jacobi polynomials and their special and degenerate cases. The operator~$L$ is up to a scaling factor, the basic hypergeometric operator;

\item[(ii)] $P_n(x)$ are the ordinary Jacobi polynomials and their degenerate cases. The operator~$L$ coincides with the ordinary hypergeometric operator;

\item[(iii)] $P_n(x)$ are the little~$-1$ Jacobi polynomials. The operator~$L$ coincides with the Dunkl type dif\/ferential operator of the f\/irst order~\cite{VZ-little}.
\end{enumerate}

In the present paper we consider the general Newtonian basis and derive the necessary and suf\/f\/icient conditions that the coef\/f\/icients~$\lambda_n$, $\tau_n$ and $a_n$ must satisfy for the associated polynomials to be orthogonal. Provided the eigenvalues $\lambda_n$ are taken to be of the ``classical'' form, we also show that this approach leads to all the polynomials of the Askey scheme (including their $q = -1$ limits). By ``classical'' we mean the Askey--Wilson grid $\lambda_n = C_1 q^n + C_2 q^{-n}+C_0$ and their degenerate forms: the quadratic grid $\lambda_n = C_1 n^2 + C_2 n + C_0$ and the Bannai--Ito grid $\lambda_n = (-1)^n (C_1 n+ C_2) + C_0$ with arbitrary constants $C_0$, $C_1$, $C_2$. We propose also a~new classif\/ication scheme of the Askey tableau elements based only on the knowledge of the coef\/f\/icients~$\lambda_n$ and~$a_n$. This classif\/ication method is purely algebraic and does not depend on any particular choice of the operator~$L$ (dif\/ferential or dif\/ference). This provides a natural ground for the hierarchy of the polynomials in the Askey tableau.

 Note that the study of orthogonal polynomials $P_n(x)$ having an expansion~\eqref{P-hyp-phi} with respect to the Newtonian basis~$\varphi_n(x)$ was initiated in the pioneering paper of Geronimus~\cite{Ger2}. For ortho\-go\-nal polynomials of several variables, the Newtonian interpolation scheme plays an important role as well~\cite{Ok,Rains}.

We here focus exclusively on the inf\/inite-dimensional situation. It should be said that the f\/inite-dimensional problem has already been considered in \cite{Ter,Ter3,Ter5,Ter4,Ter2}.
Indeed for f\/inite~$N$, the polynomials of hypergeometric type with respect to a~Newtonian basis are displayed in~\cite{PTN}.

 Note also that in the special case when all interpolation nodes $a_k$ are distinct $a_i \ne a_k$, $i \ne k$ the abstract hypergeometric polynomials have a~specif\/ic duality property. Indeed, let us def\/ine the polynomials $P_n^*(x)$, $n=0,1,\dots$ by the formula
\begin{gather*}
P_n^*(x) = \sum_{s=0}^n \frac{(a_n - a_0)(a_n-a_1)\cdots (a_n - a_{s-1}) \tilde \phi_s(x)}{\tau_1 \tau_2 \cdots \tau_s}, 
\end{gather*}
where $\tilde \phi_n(x)$ is the (dual) Newtonian basis def\/ined as
\begin{gather*}
\tilde \phi_n(x) = (x-\lambda_0)(x-\lambda_1) \cdots (x-\lambda_{n-1}). 
\end{gather*}
Then from~\eqref{P-hyp-phi} and \eqref{expr_W} we have (assuming that $W_{n0}=1$)
\begin{gather}
P_n(a_k) = P_k^*(\lambda_n). \label{LD}
\end{gather}
Property \eqref{LD} corresponds to the duality property of polynomials proposed by Leonard~\cite{Leonard}.

 In the case when both polynomial systems $P_n(x)$ and $P_n^*(x)$ are assumed to be orthogonal, it is possible to give a full classif\/ication of all such families of polynomials. For the f\/inite-dimensional case this was done by Leonard in~\cite{Leonard}; Bannai and Ito in~\cite{BI} then extended this result by including the inf\/inite-dimensional case.

The paper is organized as follows. In Section~\ref{section2} we derive the necessary and suf\/f\/icient conditions for the polynomials~$P_n(x)$ to be orthogonal. These conditions are presented in the form of a system of equations for the matrix~$Q_n^{(k)}$ of reduced moments.

In Section~\ref{section3} we show that for the ``classical'' expressions of the eigenvalues~$\lambda_n$, the remaining coef\/f\/icients~$\tau_n$ and~$a_n$ can be found explicitly. This leads to all entries of the Askey tableau.

In Section~\ref{section4} we propose a new classif\/ication scheme of the polynomials from the Askey tableau. This classif\/ication scheme dif\/fers from the known ones and is based only on the knowled\-ge of the explicit expression of the coef\/f\/icients~$\lambda_n$ and~$a_n$.

\section{Necessary and suf\/f\/icient conditions for orthogonality}\label{section2}

One way to derive the necessary and suf\/f\/icient conditions for the polynomials $P_n(x)$ to be orthogonal is to exploit the three-term recurrence relation~\eqref{rec_P}.

Assume that the polynomials~$P_n(x)$ are orthogonal which means that these polynomials satisfy~\eqref{rec_P}.

Using formulas \eqref{P-hyp-phi} and~\eqref{rec_P}, we obtain a system of equations for the coef\/f\/icients $W_{nk}$:
\begin{gather}
W_{n+1,k} + (b_n-a_{k}) W_{n,k} + u_n W_{n-1,k} = W_{n,k-1}, \qquad k=0,2,\dots, n+1. \label{eqs_W}
\end{gather}
For $k=n+1$, condition \eqref{eqs_W} becomes trivial: $1=1$. For $k=n$ and $k=n-1$, we obtain the following explicit expressions for the recurrence coef\/f\/icients $u_n$, $b_n$ in terms of the coef\/f\/i\-cients~$\lambda_n$, $\tau_n$, $a_n$
\begin{subequations}\label{ub_expl}
\begin{gather}
b_n = a_n + \frac{\tau_n}{\lambda_n - \lambda_{n-1}} - \frac{\tau_{n+1}}{\lambda_{n+1} - \lambda_{n}}, \label{b_n}\\
u_n = \frac{\tau_n(a_{n-1}-b_n)}{\lambda_n-\lambda_{n-1}} + \frac{\tau_n \tau_{n-1}}{(\lambda_n-\lambda_{n-1})(\lambda_n-\lambda_{n-2})} - \frac{\tau_n \tau_{n+1}}{(\lambda_{n+1}-\lambda_{n})(\lambda_{n+1}-\lambda_{n-1})}, \label{u_n}
\end{gather}
\end{subequations}
where we have used formulas~\eqref{W_monic}.

From \eqref{u_n} it is seen that $u_{N+1}=0$ if $\tau_{N+1}=0$. This leads to the degeneration of the corresponding orthogonal polynomials. This justif\/ies our condition $\tau_n \ne 0$ for all $n=1,2,\dots$. Nevertheless when $\tau_n \ne 0$ for $n=1,2,\dots, N$ but $\tau_{N+1}=0$ we have the special case of a f\/inite system of polynomials orthogonal on the set of points $a_0, a_1, \dots, a_N$ (in this case we should assume that all points~$a_i$ are distinct, of course). Indeed, the condition $\tau_{N+1}=0$ means that
\begin{gather*}
P_{N+1}(x) = \varphi_{N+1}(x) = (x-a_0) (x-a_1) \cdots (x-a_N). 
\end{gather*}
From the standard theory of orthogonal polynomials \cite{Chi} we have that the polynomials $P_n(x)$, $n=0,1,\dots,N$ satisfy the orthogonality relations
\begin{gather*}
\sum_{s=0}^N P_n(a_s) P_m(a_s) w_s = h_n \delta_{nm}, 
\end{gather*}
where the weights $w_s$ are expressed as
\begin{gather*}
w_s = \frac{u_1 u_2 \cdots u_N}{P'_{N+1}(a_s)P_N(a_s)}. 
\end{gather*}
Thus the condition $\tau_{N+1}=0$ leads to the special case of a f\/inite set of polynomials $P_n(x)$ orthogonal on a set $a_0,a_1, \dots, a_N$ of distinct points of the real axis. We will discuss this special case in Section~\ref{section5}.

The equations \eqref{eqs_W} corresponding to $k=n-2,n-3,\dots,0$ give an ensemble of restrictions upon the coef\/f\/icients $\lambda_n$, $\tau_n$, $a_n$. Instead of solving this system of restrictions, we shall follow a method which was successfully applied in~\cite{Zhe-hyp} to f\/ind necessary and suf\/f\/icient conditions for the orthogonality of the polynomials~$P_n(x)$.

It was shown in \cite{Zhe-hyp} that the polynomials $P_n(x)$ corresponding to the operator~$L$
are ortho\-go\-nal if and only if the operator~$L$ is symmetric on the space of polynomials. In more details, this means that for any two polynomials~$f(x)$ and $g(x)$ the condition
\begin{gather}
\langle \sigma, f(x) L g(x) \rangle = \langle \sigma, g(x) L f(x) \rangle \label{Lfg}
\end{gather}
must hold.

 Note that when $L$ is a dif\/ferential operator, condition~\eqref{Lfg} is well known~\cite{EK}. When $L$ is a~higher-order dif\/ference operator, condition~\eqref{Lfg} was derived in~\cite{Duran}. Terwilliger in~\cite{Ter2} considers this condition from the algebraic point of view in the f\/inite-dimensional case.

Condition \eqref{Lfg} is equivalent to the set of conditions
\begin{gather}
\langle \sigma, \varphi_k(x) L \varphi_n(x) \rangle = \langle \sigma, \varphi_n(x) L \varphi_k(x) \rangle, \qquad n,k =0,1,2,\dots. \label{L_nk}
\end{gather}
Taking into account relation \eqref{hyp-2-diag}, we can present~\eqref{L_nk} in the form
\begin{gather}
(\lambda_n - \lambda_k) \langle \sigma, \varphi_k \varphi_n \rangle + \tau_n \langle \sigma, \varphi_k \varphi_{n-1} \rangle - \tau_k \langle \sigma, \varphi_{k-1} \varphi_n \rangle =0 . \label{L_nk_1}
\end{gather}
We shall take the relations \eqref{L_nk_1} as the set of necessary and suf\/f\/icient conditions for the polynomials~$P_n(x)$ to be orthogonal.

It is convenient to introduce the generalized moments $\psi_n^{(k)}$
\begin{gather}
\psi_n^{(k)} \equiv \langle \sigma, \varphi_k \varphi_n \rangle, \qquad n,k=0,1,2,\dots. \label{psi_nk}
\end{gather}
By def\/inition, these moments possess the symmetry property
\begin{gather}
\psi_n^{(k)} = \psi_k^{(n)}, \qquad n,k=0,1,2,\dots . \label{sym_psi}
\end{gather}

The orthogonal polynomials $P_n(x)$ can then be presented in determinant form in terms of the generalized moments
\begin{gather*}
P_n(x;t)=\frac{1}{H_n(t)}\left |
\begin{matrix} \psi_0^{(0)} & \psi_1^{(0)} & \dots &
\psi_n^{(0)}\\ \psi_0^{(1)} & \psi_1^{(1)} & \dots & \psi_n^{(1)} \\ \dots & \dots & \dots & \dots\\
\psi_0^{(n-1)} & \psi_1^{(n-1)} & \dots & \psi_n^{(n-1)} \\ \varphi_0(x) & \varphi_1(x) & \dots & \varphi_n(x)
\end{matrix} \right |, 
\end{gather*}
where
\begin{gather*}
H_0=1, \qquad H_n=\left |
\begin{matrix} \psi_0^{(0)} & \psi_1^{(0)} & \dots &
\psi_{n-1}^{(0)}\\ \psi_0^{(1)} & \psi_1^{(1)} & \dots & \psi_{n-1}^{(1)} \\ \dots & \dots & \dots & \dots\\
\psi_0^{(n-1)} & \psi_1^{(n-1)} & \dots & \psi_{n-1}^{(n-1)}
\end{matrix} \right |, \qquad n=1,2,\dots. 
\end{gather*}
Note that the nondegeneracy condition of the linear functional $\sigma$ is equivalent to the condition
\begin{gather*}
H_n \ne 0, \qquad n=0,1,\dots . 
\end{gather*}

The relations \eqref{L_nk_1} can be presented in the form
\begin{gather}
(\lambda_n-\lambda_k) \psi_n^{(k)} + \tau_n \psi_{n-1}^{(k)} -\tau_k \psi_{n}^{(k-1)} =0. \label{L_nk_2}
\end{gather}
Moreover, from the obvious identity
\begin{gather*}
\varphi_{k+1}(x) \varphi_n(x) = \varphi_k(x) \varphi_{n+1}(x) + (a_{n}-a_{k}) \varphi_k(x) \varphi_n(x), 
\end{gather*}
we obtain the relation
\begin{gather}
\psi_n^{(k+1)} = \psi_{n+1}^{(k)} +(a_{n}-a_{k}) \psi_n^{(k)}. \label{psi_psi}
\end{gather}
Consider relation \eqref{L_nk_2} for $k=0$. Due to the conditions $\lambda_0=\tau_0=0$ we have
\begin{gather}
\lambda_n \psi_n^{(0)} = - \tau_n \psi_{n-1}^{(0)}. \label{psi_psi_0}
\end{gather}
We already know that $\lambda_n \ne 0$ for $n=1,2,\dots$. We now show that similarly $\tau_m \ne 0$ for $n=1,2,\dots$. Indeed, assume that~$\tau_j=0$ for some positive integer~$j$. Then from~\eqref{psi_psi_0} we obtain that necessarily $\psi_n^{(0)}=0$ for $n=j, j+1, j+2, \dots$. From \eqref{psi_psi} we have that $\psi_n^{(k)}=0$ for $n \ge j$ and for all $k=0,1,2,\dots$. But then $H_n=0$ for $n \ge j+1$ which would imply that the functional $\sigma$ is degenerate. We thus see that $\tau_n \ne 0$ for $n>0$. Moreover, from the same relation~\eqref{psi_psi_0}, it is seen that $\psi_0^{(0)} \ne 0$, $n=0,1,2,\dots$.

We can put $\psi_0^{(0)}=1$ (this is merely the standard normalization condition for the functional~$\sigma$). Then from~\eqref{psi_psi_0} we f\/ind the explicit expression for $\psi_n^{(0)}$:
\begin{gather}
\psi_n^{(0)} = (-1)^n \frac{\tau_1 \tau_2\cdots \tau_n}{\lambda_1 \lambda_2 \cdots \lambda_n}. \label{psi0_exp}
\end{gather}
From relation \eqref{psi_psi}, step-by-step, we can f\/ind the explicit expressions for all the generalized moments $\psi_n^{(1)}, \psi_n^{(2)}, \dots, \psi_n^{(k)}, \dots$. For example
\begin{gather*}
\psi_n^{(1)} = \left(a_{n} -\frac{\tau_{n+1}}{\lambda_{n+1}} \right) \psi_n^{(0)}. 
\end{gather*}
It is convenient to present $\psi_n^{(k)}$ as follows
\begin{gather}
\psi_n^{(k)}= Q_n^{(k)} \psi_n^{(0)} \psi_0^{(k)}, \label{psi_Q}
\end{gather}
where $Q_n^{(0)}=Q_0^{(k)}=1$.

Obviously, the reduced moments $Q_n^{(k)}$ are symmetric with respect to the variables~$n$,~$k$
\begin{gather*}
Q_n^{(k)} = Q_k^{(n)}. 
\end{gather*}
Moreover, relations \eqref{psi_psi} and \eqref{L_nk_2} become
\begin{subequations}\label{eq_Q}
\begin{gather}
\zeta_{k}Q_n^{(k+1)} = \zeta_{n}Q_{n+1}^{(k)} + (a_{k}-a_{n}) Q_n^{(k)} \label{eq_Q_1}\\
\lambda_n(Q_n^{(k)}-Q_{n-1}^{(k)}) = \lambda_k (Q_{n}^{(k)}-Q_{n}^{(k-1)}), \label{eq_Q_2}
\end{gather}
\end{subequations}
where
\begin{gather}
\zeta_n = \frac{\tau_{n+1}}{\lambda_{n+1}}, \qquad n=0,1,2,\dots. \label{zeta} \end{gather}
Relation~\eqref{eq_Q_1} allows to determine unambiguously the quantities $Q_n^{(k)}$, $k=1,2,\dots$ iteratively from the initial value $Q_n^{(0)}=0$.

The f\/irst two quantities are
\begin{gather*}
Q_n^{(1)}= \frac{\zeta_n-a_n}{\zeta_0} 
\end{gather*}
and
\begin{gather*}
Q_n^{(2)}= \frac{(a_1-a_n)(\zeta_n-a_n)}{\zeta_0 \zeta_1} + \frac{\zeta_n(\zeta_{n+1}-a_{n+1})}{\zeta_0 \zeta_1}. 
\end{gather*}
Relations~\eqref{eq_Q_2} give then additional conditions upon the quantities $Q_n^{(k)}$. These conditions are equivalent to the set of necessary and suf\/f\/icient conditions~\eqref{L_nk_1} for the orthogonality of the polynomials $P_n(x)$. Indeed, we already derived conditions \eqref{eq_Q} from the necessary and suf\/f\/icient conditions~\eqref{L_nk_2}. Conversely, assume that $Q_n^{(k)}$ is a solution of conditions~\eqref{eq_Q} for $k,n=0,1,2,\dots$ with the initial conditions $Q_n^{(0)}=Q_0^{(k)}=1$ and with the symmetry condition $Q_n^{(k)} = Q_k^{(n)}$. Then one can construct the moments~$\psi_n^{(k)}$ from~\eqref{psi_Q} and~\eqref{psi0_exp}. These moments will satisfy relations~\eqref{psi_psi} and~\eqref{sym_psi} which guarantee that the moments~$\psi_n^{(k)}$ are compatible with their def\/inition \eqref{psi_nk}. Moreover, the moments $\psi_n^{(k)}$ will satisfy the condition~\eqref{L_nk_2} for all $n,k=0,1,2,\dots$. The latter condition is equivalent to the necessary and suf\/f\/icient conditions~\eqref{L_nk}. Hence relations~\eqref{eq_Q} (together with the symmetry and initial conditions) are equivalent to the necessary and suf\/f\/icient conditions~\eqref{L_nk}.

For example, the corresponding condition coming from~\eqref{eq_Q_2} for $k=1$ looks like
\begin{gather}
\lambda_n y_{n-1} +(\lambda_1-\lambda_n) y_n + \tau_1 =0, \label{EQ_Q_1}
\end{gather}
where we have introduced the new variable
\begin{gather}
y_n = a_n -\zeta_n \label{y_def}
\end{gather}
instead of $\zeta_n$.

Similarly, for $k=2$ we obtain the relation
\begin{gather}
(y_{n+1}-y_n)(\lambda_n-\lambda_2) a_n - \lambda_n (y_n-y_{n-1})a_{n-1}\nonumber \\
\qquad {} =\lambda_n y_n (y_{n+1}-y_{n-1}) + a_1 \lambda_n (y_{n-1}-y_n) - \lambda_2 y_n (y_{n+1}-y_1). \label{EQ_k=2}
\end{gather}

It is interesting to note that the reduced moments $Q_n^{(k)}$ satisfy 3-term recurrence relations resembling the recurrence relations of the orthogonal polynomials.

Indeed, in \eqref{eq_Q_1} we can shift $k \to k-1$ and then substitute the expression for $Q_n^{(k-1)}$ given in~\eqref{eq_Q_2}. We then obtain the recurrence relation
\begin{gather}
\frac{\tau_{n+1}}{\lambda_{n+1}} (\lambda_k-\lambda_{n+1}) Q_{n+1}^{(k)} + \left(\tau_{n+1}- \tau_k +(\lambda_n-\lambda_k)(a_n-a_{k-1}) \right) Q_n^{(k)} \nonumber \\
\qquad {} +\lambda_n (a_{k-1}-a_n) Q_{n-1}^{(k)} =0, \label{rec_Q1}
\end{gather}
where the superscript $k$ is the same and the subscript $n$ takes the values $n-1$, $n$, $n+1$.

Similarly, in \eqref{eq_Q_2} we can shift $k \to k+1$ and then substitute the expression for $Q_n^{(k+1)}$. We obtain another recurrence relation:
\begin{gather}
\frac{\tau_{n+1}}{\lambda_{n+1}} (\lambda_{n}-\lambda_{k+1}) Q_{n+1}^{(k)} + \left(\tau_{k+1}- \tau_n +(\lambda_n-\lambda_{k+1})(a_k-a_n) \right) Q_n^{(k)} \nonumber \\
\qquad {} +\lambda_n (a_{n-1}-a_{k}) Q_{n-1}^{(k)} =0. \label{rec_Q2}
\end{gather}
Relations \eqref{rec_Q1} and \eqref{rec_Q2} are not independent due to relations~\eqref{eq_Q}.

\section{General solution for prescribed classical eigenvalues}\label{section3}

\subsection[The case of linear $\lambda_n$]{The case of linear $\boldsymbol{\lambda_n}$}\label{section3.1}

The general solution of the necessary and suf\/f\/icient conditions~\eqref{eq_Q} is rather complicated. Instead we present here a special solution of this problem, starting with prescribed dependence of $\lambda_n$ on $n$. Let us consider f\/irst the simplest case of linear dependence: $\lambda_n = \alpha n + \beta$. Clearly, the scaling transformation $\lambda_n \to \kappa \lambda_n$, $\tau_n \to \kappa \tau_n$ leaves the problem invariant, so we can choose~$\alpha=1$. Moreover, the initial condition $\lambda_0=0$ leads to $\beta=0$. Hence we may take~$\lambda_n=n$ without loss of generality.

Consider equation \eqref{EQ_Q_1} with $\lambda_n=n$.
The general solution of this equation is elementary
\begin{gather}
y_n = \gamma n - \tau_1, \label{y_tau_lin} \end{gather}
where $\gamma$ is an arbitrary constant.

Substituting expression \eqref{y_tau_lin} into equation \eqref{EQ_k=2} we obtain the equation for the grid $a_n$
\begin{gather*}
\frac{a_{n-1}}{(n-1)(n-2)} - \frac{a_{n}}{n(n-1)} = \frac{a_{1}}{(n-1)(n-2)}. 
\end{gather*}
The solution of this equation is immediate:
\begin{gather}
a_n = a_1 n + \alpha n(n-1) \label{gen_sol_a_lin} \end{gather}
with an arbitrary constant~$\alpha$.

From \eqref{y_def} we obtain the solution for $\tau_n$:
\begin{gather}
\tau_{n+1} = (n+1) \left(\tau_1 +(a_1-\gamma) n + \alpha  n (n-1) \right). \label{gen_sol_tau_lin}
\end{gather}
We have thus found the general solution of the problem for the special case $\lambda_n=n$. Instead of checking the compatibility of all conditions~\eqref{eq_Q} with the obtained solution it is suf\/f\/icient to notice that this solution describes orthogonal polynomials that are already known.

Indeed, if $\alpha=a_1=0$, then $a_n=0$, $n=0,1,2,\dots$ and we obtain the monomial basis $\phi_n(x) = x^n$. For~$\tau_n$ we have from \eqref{gen_sol_tau_lin}:
\begin{gather*}
\tau_n = n \left(\tau_1 -\gamma (n-1) \right). 
\end{gather*}
The corresponding operator $L$ is the second-order dif\/ferential operator
\begin{gather}
L = -\gamma  x \partial_x^2 + (x-\tau_1) \partial_x . \label{Lag_L}
\end{gather}
It is well known that the operator \eqref{Lag_L} has Laguerre polynomials as eigenfunctions~\cite{KLS}.

If $\alpha=0$ but $a_1 \ne 0$, then the grid is linear:
\begin{gather*}
a_n = a_1 n. 
\end{gather*}
The coef\/f\/icient $\tau_n$ is a quadratic function in $n$:
\begin{gather*}
\tau_n = n(\tau_1 +(a_1-\gamma)(n-1)). 
\end{gather*}
This case corresponds to all classical polynomials on the uniform grid: Krawtchouk, Meixner, Meixner--Pollaczek and Charlier polynomials.

Finally, if $\alpha \ne 0$ we have the quadratic grid \eqref{gen_sol_a_lin}. This corresponds to the dual Hahn polynomials (including the continuous dual Hahn polynomials)~\cite{KLS}.

The analysis of the case where $\lambda_n$ is linear in $n$ expounds all known classical orthogonal polynomials with such a spectrum. The only exception are the Hermite polynomials, because the eigenvalue equation for these polynomials in the monomial basis~\cite{KLS} $x^n$ contains the terms~$x^n$ and~$x^{n-2}$ instead of~$x^n$ and~$x^{n-1}$ as required in our approach.

\subsection[The case of quadratic $\lambda_n$]{The case of quadratic $\boldsymbol{\lambda_n}$}\label{section3.2}

Consider the case when $\lambda_n$ is a quadratic polynomial in $n$. Because $\lambda_0=0$ it is suf\/f\/icient to take
\begin{gather}
\lambda_n = n(n+\alpha) \label{quad_lambda}
\end{gather}
with some real parameter $\alpha$. The only restriction is $\alpha \ne -1,-2, -3, \dots$ because otherwise $\lambda_j=0$ for some positive integer $j$ and that is forbidden by hypothesis since we are assuming that the spectrum is non degenerate.

Consider f\/irst equation \eqref{EQ_Q_1}. Substituting expression \eqref{quad_lambda}, we get
\begin{gather}
\frac{(n+\alpha)y_{n-1}}{n-1} - \frac{(n+\alpha+1)y_{n}}{n} + \tau_1 \left( \frac{1}{n-1} - \frac{1}{n} \right) =0 \label{eq_1_quad}
\end{gather}
(with $y_n$ again as in \eqref{y_def}). The solution of \eqref{eq_1_quad} is immediately found to be:
\begin{gather}
y_n =\frac{\gamma n - \tau_1}{n+\alpha+1}, \label{y_quad} \end{gather}
where $\gamma$ is an arbitrary constant.

Equation \eqref{EQ_k=2} now reads
\begin{gather*}
\kappa (n a_{n-1} - (n-2) a_n -a_1 n ) =0, 
\end{gather*}
where
\begin{gather*}
\kappa = (\alpha+2 ) (\tau_1 +\gamma (\alpha+1) ). 
\end{gather*}
Note that $\kappa \ne 0$. Indeed, $\alpha=-2$ means that $\lambda_2=0$ which is forbidden since the spectrum must be non degenerate. If $\tau_1 +\gamma (\alpha+1)=0$, it is easy to show from equation \eqref{eq_Q_1} that all moments are degenerate $Q_n^{(1)}=Q_n^{(2)} =\cdots = \operatorname{const}$ which is again forbidden for degenaracy reasons. We thus conclude that $\kappa \ne 0$ and hence have the equation
\begin{gather}
n a_{n-1} - (n-2) a_n -a_1 n =0 \label{eq_a_quadr}
\end{gather}
which has the general solution
\begin{gather}
a_n = a_1 n + \beta n (n-1) \label{sol_a_quad}
\end{gather}
with an arbitrary constant $\beta$. (The value of $a_1$ is arbitrary as well, it can be considered as the initial condition for~$a_n$). We thus obtained a quadratic grid for~$a_n$.

Finally, from \eqref{y_def} and \eqref{y_quad}, we f\/ind the explicit expression for $\tau_n$:
\begin{gather*}
\tau_{n+1} = \lambda_{n+1} a_{n} + (n+1) (\tau_1 - \gamma n), 
\end{gather*}
where $\lambda_n$ and $a_n$ are given by~\eqref{quad_lambda} and~\eqref{sol_a_quad}.

In general, $\tau_n$ is a polynomials of 4-th degree in $n$. In this case the grid $a_n$ is quadratic and the corresponding polynomials are the Wilson--Racah polynomials. When $\beta =0$ but $a_1 \ne 0$, then the grid $a_n$ is linear, $\tau_n$ is a polynomial of third degree and the corresponding polynomials are the Hahn polynomials. Finally, when $a_1=\beta=0$, then $a_n =0$ for all $n$. The basis $\varphi_n(x) =x^n$ is the monomial one and the corresponding polynomials $P_n(x)$ are the Jacobi polynomials.

\subsection{The Askey--Wilson grid}\label{section3.3}

Consider the eigenvalues
\begin{gather}
\lambda_n = \big(1-q^n\big)\big(\alpha-q^{-n}\big) \label{AW_l}
\end{gather}
with an arbitrary real parameter $\alpha$.
These eigenvalues correspond to the Askey--Wilson polynomials.

 In this case, the general solution to~\eqref{EQ_Q_1} for $y_n$ is{\samepage
\begin{gather}
y_n = \frac{\tau_1 (q-1)^{-1} q^{-n} + \gamma (q^{-n}-1)}{\alpha-q^{-n-1}} \label{y_AW}
\end{gather}
with $\gamma$ an arbitrary parameter.}

Substituting then expression~\eqref{y_AW} into~\eqref{EQ_k=2}, we obtain the equation for the unknown grid~$a_n$:
\begin{gather}
\big(q^{n-1}-q\big)a_n -\big(q^n-1\big)a_{n-1}=a_1\big(1-q^n\big). \label{eq_a_AW}
\end{gather}
This is a simple f\/irst order dif\/ference equation which has for general solution
\begin{gather}
a_n = \big(1-q^{-n}\big) \big(a_1q (q-1)^{-1} + \nu \big(q^{n-1}-1\big)\big) \label{a_n_AW}
\end{gather}
with $\nu$ an arbitrary parameter.

Expression \eqref{a_n_AW} corresponds to the grid of the Askey--Wilson polynomials. Indeed, one checks directly that the recurrence coef\/f\/icients $b_n$,~$u_n$ constructed with the help of formu\-las~\eqref{ub_expl} correspond to the most general Askey--Wilson polynomials. There are 4 free parameters: $\tau_1$,~$a_1$ and the integration constants~$\gamma$,~$\nu$. These parameters correspond to the 4~parameters of the Askey--Wilson polynomials.

\subsection{Bannai--Ito grid}\label{section3.4}

Finally, consider the case of the Bannai--Ito eigenvalues \cite{TVZ}
\begin{gather}
\lambda_n = (-1)^n (n+\alpha) -\alpha \label{BI_l} \end{gather}
with an arbitrary parameter $\alpha$. Note that for real $\alpha$ the even and the odd eigenvalues form two equidistant lattices each with step~2. These sublattices are separated by the gap $\Delta = 2 \alpha+1$. For the special case $\alpha=1/2$, these two sublattices can be combined into one equidistant lattice $\lambda_n = \pm 2n$, $n=0, 1, 2,\dots$.

The eigenvalues \eqref{BI_l} can be obtained from a specif\/ic limit of the AW eigenvalues \eqref{AW_l} (see, e.g.,~\cite{TVZ}). It is instructive however to start from~\eqref{BI_l} to see how the analysis proposed here applies directly to the BI lattice.

The f\/irst equation, namely, \eqref{EQ_Q_1} can be split into two equations according to the parity of~$n$
\begin{gather}
 2n y_{2n-1} -(2n+2\alpha+1) y_{2n} + \tau_1 =0, \\
 2n y_{2n+1} -(2n+2\alpha+1) y_{2n} + \tau_1 =0. 
\end{gather}
The solution of these equations is readily seen to be
\begin{gather}
y_{2n} = \frac{2 \gamma n + \tau_1}{2n + 2 \alpha+1}, \qquad y_{2n+1}=\gamma, \label{y_BI} \end{gather}
where $\gamma$ is an arbitrary parameter.

Substituting \eqref{y_BI} for $y_n$ into \eqref{EQ_k=2}, we obtain the following equation for the unknown grid~$a_n$:
\begin{gather}
\kappa \left((n-1) a_{2n} + n a_{2n-1}-a_1 n \right)=0, \qquad \kappa(a_{2n+1} + a_{2n}-a_1)=0, \label{eq_a_BI} \end{gather}
where
\begin{gather*}
\kappa = \gamma (2\alpha+1) - \tau_1. 
\end{gather*}
As in the previous cases the constant $\kappa$ should be nonzero to avoid degeneracies. The solution of the equations \eqref{eq_a_BI} is elementary
\begin{gather*}
a_{2n} = 2 \nu n, \qquad a_{2n+1} = a_1 - 2 \nu n, \qquad n=0,1,2,\dots, 
\end{gather*}
where $\nu$ is an arbitrary parameter.

The coef\/f\/icients $\tau_n$ are then expressed as
\begin{gather}
\tau_n = \lambda_n (a_{n-1}-y_{n-1}), \qquad n=1,2,\dots. \label{tau_BI} \end{gather}
Using the expressions for~$a_n$,~$\lambda_n$,~$\tau_n$, we can calculate the recurrence coef\/f\/icients $b_n$, $u_n$ with the help of formulas~\eqref{ub_expl}. These coef\/f\/icients coincide with the recurrence coef\/f\/icients of the Bannai--Ito polynomials~\cite{TVZ} for an arbitrary choice of the parameters~$a_1$, $\tau_1$,~$\gamma$,~$\nu$.

 We have thus considered all possible choices for the classical spectrum $\lambda_n$; in all cases the orthogonal polynomials that are found coincide with the corresponding polynomials from the Askey tableau. It is quite natural to conjecture that there are no other (``non-classical'') expressions for~$\lambda_n$ which give rise to orthogonal polyno\-mials. This means that only the polynomials from the Askey tableau do satisfy the conditions of the Newtonian ``hypergeometric'' polyno\-mials. This conjecture remains however an open problem.

\section{Algebraic classif\/ication of the elements of the Askey tableau}\label{section4}

All the elements of the Askey tableau can be characterized naturally within the unif\/ied algebraic ``abstract'' hypergeometric scheme of\/fered here. We observed that the broad features of the spectrum $\lambda_n$ provide a f\/irst categorization. Indeed, we already know that all classical polynomials from the ``abstract'' hypergeometric scheme can be classif\/ied by the eigenvalues $\lambda_n$. There are 3~main classes: one corresponding to a~generic~$q$ (Askey--Wilson class) and two associated to the special cases, $q=1$ (Racah--Wilson class) and $q=-1$ (Bannai--Ito class). Within each of these classes, the interpolation grids $a_n$ further split into 3 dif\/ferent types. This gives 9 basic classes of classical hypergeometric polynomials.

Consider for example the general Askey--Wilson class. The grid $a_n$ may be taken to be the most general~-- this gives the generic AW polynomials. If the grid is exponential (this corresponds to $\nu=0$), we obtain the big $q$-Jacobi polynomials. Finally, if the grid degenerates to a single point (i.e. when $\nu=a_1=0$), we obtain the little $q$-Jacobi polynomials.

For the Racah--Wilson class ($q=1$), we have 3 subclasses according to the particular expression \eqref{eq_a_quadr} of the interpolation grid $a_n$. If both parameters $a_1$ and $\beta$ are nonzero, we are then dealing with the general Racah--Wilson polynomials. If $\beta=0$ but $a_1 \ne 0$, the interpolation grid becomes linear and we have the Hahn polynomials (including the case of the continuous Hahn polynomials). Finally, if $a_1=\beta=0$, the grid is then completely degenerate $a_n=0$, the corresponding basis is the monomial one $\phi_n(x) =x^n$ and we have the Jacobi polynomials.

For the Bannai--Ito class we have again 3 subcategories depending on the expression of $a_n$. If both parameters $a_1$ and $\nu$ are nonzero, then the grid $a_n$ is the general Bannai--Ito grid and this leads to the Bannai--Ito polynomials. If $\nu=0$ but $a_1 \ne 0$, we have a degenerate grid which consists of only two distinct points: $a_{2n}=0$ and $a_{2n+1}=a_1$; this case corresponds to the big~$-1$ Jacobi polynomials. Finally, if $a_1=\nu=0$ then the grid becomes completely degenerate $a_n=0$ and the corresponding polynomials are the little~$-1$ Jacobi polynomials.

One can consider further subdivisions inside these 9~classes. For example, in the Racah--Wilson class one may consider the case of linear eigenvalues $\lambda_n$; one then obtains the classical orthogonal polynomials of the Hahn class and their degenerate types (Meixner, Krawtchouk, Charlier and Laguerre). Similarly, if one takes the case of the exponential AW-grid corresponding to $\alpha=0$ in~\eqref{AW_l}, one is led to the $q$-Hahn polynomials and their degenerate and special cases.

\section{Conclusions}\label{section5}

In summary, we derived the necessary and suf\/f\/icient conditions for the orthogonality of the ``hypergeometric'' polynomials with respect to the Newton interpolation basis. These conditions are presented in the form of a system of equations~\eqref{eq_Q} involving the reduced moments~$Q_n^{(k)}$. When the eigenvalues~$\lambda_n$ are given by ``classical'' expressions (i.e., trigonometric, quadratic or Bannai--Ito grid), we have showed that the corresponding polynomials~$P_n(x)$ coincide with the entries of the Askey tableau. This gives a simple and natural new classif\/ication of all polynomials of the Askey tableau.

Let us remark that the monomial basis (which generates Jacobi, Laguerre and Hermite polynomials) is a degeneration of the Newton interpolation basis. As pointed out by a~referee, let us mention that the little $q$-Jacobi polynomials can also be viewed as a special case of the big $q$-Jacobi polynomials (see, e.g., formula~(III.7) in~\cite{GR}) and hence also associated with to the same Newtonian basis as the latter.

An unresolved issue is to prove that the general problem with all unknown coef\/f\/i\-cients $\lambda$,~$\tau_n$,~$a_n$ gives in all likelihood the same solution (i.e., the polynomials of the Askey tableau). Another interesting question is to f\/ind the explicit expression of the moments~$Q_n^{(k)}$ for all polynomial families of the Askey tableau. Indeed, formulas~\eqref{rec_Q1} and~\eqref{rec_Q2} can be considered as three-term recurrence relations for these moments; one could expect that they are related to the recurrence relations of the corresponding orthogonal polynomials.

The problem considered here admits several possible generalizations. We mention only one of them. Fix some positive integer $j$ and assume that the operator $L$ satisfy the following more general property in the Newtonian basis~$\varphi_n(x)$:
\begin{gather}
L \varphi_n(x) = \sum_{s=0}^j \tau_{ns} \varphi_{n-s}(x) \label{L_gen_j}
\end{gather}
with some coef\/f\/icients $\tau_{ns}$ such that $\tau_{n0}=\lambda_n$. Then the eigenvalue problem~\eqref{LPP} generates again a set of unique monic polynomials $P_n(x)$. The issues examined in the present paper correspond to the simplest nontrivial case $j=1$. For $j=2,3,\dots$ we can consider similarly the orthogonality of the corresponding polynomials $P_n(x)$. Of course, all ``classical'' orthogonal polynomials are obtained as special cases for every $j$. Indeed, let $L$ be the operator associated to the classical polynomials $P_n(x)$, then the operator $\tilde L = \xi_j L^j + \xi_{j-1} L^{j-1} + \cdots +\xi_1 L + \xi_0$ (i.e., any polynomial of the operator~$L$ of degree~$j$) will trivially have the same orthogonal polynomials as solutions of the eigenvalue problem
\begin{gather*}
\tilde L P_n(x) = \tilde \lambda_n P_n(x). 
\end{gather*}
The nontrivial question is the existence of non-classical orthogonal polynomials satisfying~\eqref{LPP} with ``higher-order'' operators~\eqref{L_gen_j}. Such non-classical polynomials appear, e.g., in the solution of the Krall problem that describes all the orthogonal polynomials satisfying higher-order dif\/ferential equations~\cite{EK}. We can therefore expect the appearance of generalized Krall polynomials (e.g., those introduced by Dur\'an~\cite{Duran}) as solutions of this problem.

\subsection*{Acknowledgments}

The authors thank the referees and the editors for valuable remarks and suggestions. AZ thanks to the Centre de Recherches Math\'ematiques (Universit\'e de Montr\'eal) for hospitality. The authors would like to thank S.~Tsujimoto for stimulating
discussions. The research of LV is supported in part by a research grant from the Natural Sciences and Engineering Research Council
(NSERC) of Canada.

\pdfbookmark[1]{References}{ref}
\LastPageEnding

\end{document}